\documentclass[twocolumn]{autart}    

\usepackage{amsmath,amssymb,epsf,epsfig,times}
\usepackage[all]{xy}
\usepackage{graphicx,color}
\usepackage{subfigure}
\usepackage{url}
\usepackage{cite}
\usepackage{epstopdf}

\newtheorem{theorem}{Theorem}[section]
\newtheorem{lemma}{Lemma}[section]

\def\QED{\mbox{\rule[0pt]{1.5ex}{1.5ex}}}
\def\endproof{\hspace*{\fill}~\QED\par\endtrivlist\unskip}
\newcommand{\re}{\mathbb{R}}

\newcommand{\pde}[2]{\frac{\partial #1}{\partial #2}}

\newcommand{\defeq}{\stackrel{\triangle}{=}}

\newtheorem{assumption}[theorem]{Assumption}

\newtheorem{remark}[theorem]{Remark}

\newcommand{\OMIT}[1]{}

\begin{document}

\begin{frontmatter}

\title{Fully Distributed Flocking with a Moving Leader for Lagrange Networks with Parametric Uncertainties \thanksref{footnoteinfo}} 

\thanks[footnoteinfo]{
A preliminary version of the paper has appeared at the ACC 2014.}

\author[Paestum]{Sheida Ghapani}\ead{sghap001@ucr.edu}~,
\author[Paestum1]{Jie Mei}\ead{meij.hit@gmail.com}~,
\author[Paestum]{Wei Ren}\ead{ren@ee.ucr.edu}~,
and  
\author[Paestum2]{Yongduan Song}\ead{ydsong@cqu.edu.cn}

\address[Paestum]{Department of Electrical and Computer Engineering, University of California, Riverside, USA}  
\address[Paestum1]{School of Mechanical Engineering and Automation, Harbin Institute of Technology Shenzhen Graduate School, Shenzhen, China}             
\address[Paestum2]{School of Automation, Chongqing University, Chongqing, China}             



\begin{abstract}                          
This paper addresses the leader-follower flocking problem with a moving leader for networked Lagrange systems with parametric uncertainties under a proximity graph. Here a group of followers move cohesively with the moving leader to maintain connectivity and avoid collisions for all time and also eventually achieve velocity matching. In the proximity graph, the neighbor relationship is defined according to the relative distance between each pair of agents. Each follower is able to obtain information from only the neighbors in its proximity, involving only local interaction. We consider two cases: i) the leader moves with a constant velocity, and ii) the leader moves with a varying velocity. In the first case, a distributed continuous adaptive control algorithm accounting for unknown parameters is proposed in combination with a distributed continuous estimator for each follower. In the second case, a distributed discontinuous adaptive control algorithm and estimator are proposed. Then the algorithm is extended to be fully distributed with the introduction of gain adaptation laws. In all proposed algorithms, only one-hop neighbors' information (e.g., the relative position and velocity measurements between
the neighbors and the absolute position and velocity measurements) is required, and flocking is achieved as long as the connectivity and collision avoidance are ensured at the initial time and the control gains are designed properly. Numerical simulations are presented to illustrate the theoretical results.
\end{abstract}

\end{frontmatter}

\section{Introduction} \label{sec:Introduction}

A multi-agent system is defined as a collection of autonomous agents which are able to interact with each other or with their environments to solve problems that are difficult or impossible for an individual agent. In a multi-agent system, the agents often act in a distributed manner to complete global tasks cooperatively with only local information from their neighbors so as to increase flexibility and robustness.

The collective behavior can be observed in nature like flock of birds, swarm of insects, and school of fish. In \cite{reynolds1987flocks}, three heuristic rules are characterized for the flocking of multi-agent systems, namely, flock centering, collision avoidance and velocity matching. In \cite{tanner2007flocking}, a flocking algorithm is introduced for a group of agents when there is no leader. A theoretical framework is proposed in \cite{olfati2006flocking} to address the flocking problem with a leader, which has a constant velocity and is a neighbor of all followers.
Ref. \cite{su2009flocking} considers both cases where the leader has a constant and a varying velocity. When the leader has a constant velocity, \cite{su2009flocking} relaxes the constraint that the leader is a neighbor of all followers. However, in the case where the leader has a varying velocity, it still requires that the leader be a neighbor of all followers. Unfortunately, this is an unrealistic restriction on the distributed control design, especially when the number of the followers becomes large.
In \cite{cao2012distributed}, distributed control algorithms for swarm tracking are studied via a variable structure approach, where the moving leader is a neighbor of only a subset of the followers. In \cite{li2013flocking}, the flocking control and communication optimization problem is considered for multi-agent systems in a realistic communication environment and the desired separation distances between neighboring agents is calculated in real time.

Note that all above references focus on linear multi-agent systems with single- or double-integrator dynamics. However, in reality, many physical systems are inherently nonlinear and cannot be described by linear equations. Among the nonlinear systems, Lagrange models can be used to describe a large class of physical systems of practical interests such as autonomous vehicles, walking robots, and rotation and translation of spacecraft formation flying. But due to the existence of nonlinear terms with parametric uncertainties, the algorithms for linear models cannot be directly used to solve the coordination problem for multi-agent systems with Lagrange dynamics.

Recent results on distributed coordination of networked Lagrange systems focus on the consensus without a leader \cite{ChopraSpong06,ren2009distributed,HouChengTan09,MinSWL11,Wang13,Wang14_TAC}, coordinated tracking with one leader \cite{chung2009cooperative,mei2011distributed,dong2011consensus}, containment control with multiple leaders \cite{meng2010distributed,mei2012distributed,MeiRenChenMa13_automatica}, and flocking or swarming without or with a leader \cite{chopra2008synchronization,cheah2009region,meng2012leader}.
Ref. \cite{chopra2008synchronization} proposes a control algorithm based on potential functions for networked Lagrange systems to achieve collision avoidance and velocity matching simultaneously in both time-delay and switching-topology scenarios. However, parametric uncertainties are not considered and there is no leader. 
Ref. \cite{cheah2009region} presents a region-based shape controller for a swarm of Lagrange systems.
By utilizing potential functions, the authors design a control scheme that can force multiple robots to move as a group inside a desired region with a common velocity while maintaining a minimum distance among themselves. However, the algorithm relies on the strict assumption that all followers have access to the information of the desired region and the common velocity. A leader-follower swarm tracking framework is established in \cite{meng2012leader} in the presence of multiple leaders. However, only a compromised result can be obtained when the group dispersion, cohesion, and containment objectives are considered together. In the proposed algorithms, the variables of the estimators must be communicated among the followers. Furthermore, more information is used in the controller design, for example, the second-order derivatives of the potential functions.

In this paper we focus on the distributed leader-follower flocking problem with a moving leader for networked Lagrange systems with unknown parameters under a proximity graph defined according to the relative distance between each pair of agents. Here a group of followers move cohesively with the moving leader to maintain connectivity and avoid collisions for all time and also eventually achieve velocity matching. The leader can be a physical or virtual vehicle, which encapsulates the group trajectory. We consider two cases: i) the leader moves with a constant velocity, and ii) the leader moves with a varying velocity. In the first case, a distributed continuous adaptive control algorithm accounting for unknown parameters and a distributed continuous estimator is proposed for each follower. In the second case, we first propose a distributed discontinuous adaptive control algorithm and estimator, where we use a common control gain that is sufficiently large for all followers. Hence the system is not completely distributed. We then improve the algorithm by further proposing gain adaption schemes to implement a fully distributed algorithm.
In all proposed algorithms, only one-hop neighbors' information is used, and flocking is achieved as long as the connectivity and collision avoidance are ensured at the initial time and the control gains are designed properly.
Compared with the results in the existing literature, this paper has the following novel features.
\begin{itemize}
  \item[1)] This paper considers each agent as a nonlinear Euler-Lagrange system with parametric uncertainties and is more realistic. While in \cite{olfati2006flocking,tanner2007flocking,cao2012distributed,li2013flocking}, the agents' dynamics are assumed to be single or double integrators. The results for single- or double-integrator dynamics are not applicable to Lagrange systems with parametric uncertainties.
  \item[2)] This paper considers the combination of flocking (considering connectivity maintenance, collision avoidance, and velocity
matching with a moving leader in the meantime) and the constraint that the leader's information is available to only the followers in its proximity. The above constraint introduces further complexities since not all followers know the leader's velocity. Even for the case with single- or double-integrator agents, the problem is very challenging \cite{cao2012distributed}, not to mention the case of nonlinear Lagrange systems with parametric uncertainties. In contrast, in \cite{chopra2008synchronization}, parametric uncertainties
are not considered and there is no leader and in \cite{cheah2009region}, it is assumed that the leader's information is available to all followers (against the local interaction nature of the problem).
  \item[3)] To overcome the coexistence and coupling of the above mentioned challenges, in the current paper, we propose an adaptive control law in combination with
a new distributed estimator for each follower. The novelty of the estimators is that the partial derivatives of the potential
functions are integrated into the estimators. In \cite{cao2012distributed,meng2012leader}, the variables of the estimators must be communicated between the neighbors.
For the case of a moving leader with varying velocity, the proposed algorithms in \cite{cao2012distributed,mei2011distributed} require both one-hop and two-hop neighbors' information. In contrast, in our proposed algorithms, only one-hop neighbors' information (e.g., the relative position and velocity measurements between
the neighbors and the absolute position and velocity measurements) is required. These measurements can be obtained by the sensing devices
carried by the agents and hence the need for communication can be removed.
 Further, a fully distributed algorithm without global information is proposed in the current paper, while the results in \cite{cao2012distributed,mei2011distributed,meng2012leader} rely on some global information.
\end{itemize}

{\it Notations:} Let $\mathbf{1}_n$
denote the $n \times 1$ column vector of all ones.
Let $\lambda_{\min} (.)$ denote the minimum eigenvalue of a square real matrix with real eigenvalues. Let $\mbox{diag}(z_1,\ldots,z_p)$ be the diagonal matrix with diagonal entries $z_1$ to $z_p$. For symmetric square real matrices $A$ and $B$ with the same order, $A>B$ or equivalently $B<A$ (respectively, $A \geq B$ or equivalently $B \leq A$) means that $A-B$ is symmetric positive definite (respectively, semi-definite). Throughout the paper, we use $||\cdot||$ to denote the Euclidean norm, $\otimes$ to denote the Kronecker product, and $\mbox{sgn}(\cdot)$ to denote the $\mbox{signum}$ function defined componentwise. For a vector function $f(t):\re\mapsto\re^m$, it is said that
$f(t)\in\mathbb{L}_l$ if
$(\int_{0}^{\infty}\|f(\tau)\|^l\mbox{d}\tau)^{\frac{1}{l}}<\infty$ and
$f(t)\in\mathbb{L}_{\infty}$ if for each element of $f(t)$, noted as
$f_i(t)$, $\sup_{t \geq 0}|f_i(t)|<\infty$, $i=1,\ldots,m$.

\section{Background}

\subsection{Lagrange Dynamics}

Suppose that there exist $n+1$ agents (e.g., autonomous vehicles) consisting of one leader and $n$ followers. The leader is labeled as agent $0$ and the followers are labeled as agent $1$ to $n$. The $n$ followers are described by Lagrange equations of the form \cite{kelly2006control}
\begin{equation}\label{EL-dyn}
M_i (q_i )\ddot{q}_i+C_i (q_i,\dot{q}_i )\dot{q}_i+g_i (q_i)=u_i,  \qquad i=1,\ldots,n,
\end{equation}
where $q_i \in \mathbb{R}^p$ is the vector of generalized coordinates\footnote{In the context of autonomous vehicles, $q_i$ denotes the position of agent $i$.}, $M_i (q_i )$ is the $p \times p$ symmetric inertia matrix, $C_i (q_i,\dot{q}_i )\dot{q}_i$ is the Coriolis and centrifugal force, $g_i (q_i)$ is the vector of gravitational force, and $u_i$ is the control input. The dynamics of the Lagrange systems satisfy the following properties:

\begin{itemize}

\item[(P1)] There exist positive constants $k_{\underline{M}}, k_{\overline{M}}, k_{\overline{C}} , k_{\overline{g}}$ such that $k_{\underline{M}} I_p \le M_i (q_i ) \le k_{\overline{M}} I_p, ||C_i (q_i,\dot{q}_i )\dot{q}_i || \le k_{\overline{C}} ||\dot{q}_i ||$ and $||g_i (q_i)|| \le k_{\overline{g}}$ .
\item[(P2)] $\dot{M}_i (q_i )-2C_i (q_i,\dot{q}_i )$ is skew symmetric.
\item[(P3)] The left-hand side of the Lagrange dynamics can be parameterized, i.e., $M_i (q_i )x+C_i (q_i,\dot{q}_i )y+g_i (q_i)=Y_i (q_i,\dot{q}_i,x,y) \theta_i$, $\forall x,y \in \mathbb{R}^p$, where $Y_i \in \mathbb{R}^{p \times p_ \theta}$ is the regression matrix and $\theta_i \in \mathbb{R}^{p_\theta}$ is the unknown but constant parameter vector.

\end{itemize}

In this paper, the leader can be a physical or virtual vehicle, which encapsulates the group trajectory. The leader's position and velocity are denoted by, respectively, $q_0 \in \mathbb{R}^p$ and $\dot{q}_0 \in \mathbb{R}^p$.

\subsection{Graph Theory}

With $k$ agents in a team, a graph is used to characterize the interaction topology among the agents.
A graph is a pair $G=(V,E)$, where $V=\{1,\ldots,k\}$ is the node set and $E \subseteq V \times V$ is the edge set. In a directed graph, an edge $(j,i) \in E$ means that node $i$ can obtain information from node $j$ but not necessarily vice versa. Here node $j$ is a neighbor of node $i$. In an undirected graph $(i,j) \in E \Leftrightarrow (j,i) \in E$. A directed path in a directed graph is an ordered sequence of edges of the form $(i_1,i_2)$, $(i_2,i_3), \ldots ,$ where $i_j \in V$. A subgraph of $G$ is a graph whose node set and edge set are subsets of those of $G$.

The adjacency matrix $\mathbf{A}=[a_{ij}] \in \mathbb{R}^{k \times k}$ of the graph $G$ is defined such that the edge weight $a_{ij}=1$ if $(j,i) \in E$ and $a_{ij}=0$ otherwise. For an undirected graph, $a_{ij}=a_{ji}$. The Laplacian matrix $L=[l_{ij}] \in \mathbb{R}^{k \times k}$ associated with $\mathbf{A}$ is defined as $l_{ii}=\sum_{j \ne i} a_{ij}$ and $l_{ij}=-a_{ij}$, where $i \ne j$. For an undirected graph, $L$ is symmetric positive semi-definite \cite{chung1997spectral}.

In this paper, we assume that the neighbor relationship among the leader and the followers is based on their relative distance and hence the graph characterizing the interaction topology is a {\it proximity} graph. We also assume that the leader has no neighbor and its motion is not necessarily dependent on the followers. In particular, followers $i$ and $j$ are neighbors of each other if $||q_i-q_j||<R$ and the leader is a neighbor of follower $i$ if $||q_i-q_0||<R$, where $R$ denotes the sensing radius of the agents. Let $G_F$ be the proximity graph characterizing the interaction among the $n$ followers with the associated Laplacian matrix $L_F$. Note that by definition $G_F$ is undirected and hence $L_F$ is symmetric positive semi-definite. To
simplify our analysis, we assign an orientation to an edge by
considering one node the positive end of the edge and the other node
the negative end of the edge. We recall that the $n\times N$
incidence matrix $D_F=[d_{ik}]\in \mathbb{R}^{n\times N}$ of a graph
is defined as \cite{Biggs93}
\begin{equation}\label{D_definition}
    d_{ik}=\left\{\begin{array}{ll}
                  +1 & \mbox{if node $i$ is the positive end of the edge $\mathcal {E}_k$}, \\
                  -1 & \mbox{if node $i$ is the negative end of the edge $\mathcal {E}_k$}, \\
                  0 &  \mbox{otherwise.}
                \end{array}
                \right. \nonumber
\end{equation}
Then the Laplacian matrix of the graph can be denoted by $L_F=D_FD_F^T$.

Let $\overline{G}$ be the directed graph characterizing the interaction among the leader and the $n$ followers corresponding to $G_F$. Also let the edge weight $a_{i0}=1$ if the leader is a neighbor of follower $i$ and $a_{i0}=0$ otherwise. Define $\Lambda \defeq \mbox{diag}(a_{10},\ldots,a_{n0})$. Note that $\Lambda^2=\Lambda$ because $a_{i0}$ is either 1 or 0. Also define the {\it leader-follower topology matrix} associated with the graph $\overline{G}$ as $H=L_F+\Lambda$. It is obvious that $H$ is symmetric positive semi-definite. Before moving on, we need the following lemmas.

\begin{lemma}\cite{ren2010distributed} \label{lemma:tree}
If the leader has directed paths to all followers, the matrix $H$ is symmetric positive definite.
\end{lemma}

\begin{lemma}\label{lemma:psd}
Let $H^a$ and $H^b$ be the leader-follower topology matrix associated with, respectively the graph $\overline{G}^a$ and $\overline{G}^b$. If $\overline{G}^a$ is a subgraph of $\overline{G}^b$, then $H^a \leq H^b$.
\end{lemma}

\emph{Proof}: When $\overline{G}^a$ is a subgraph of $\overline{G}^b$, $H^b$ can be written as $H^b=H^a+P$, where $P$ is a positive semi-definite matrix. Therefore, it can be concluded that $H^a \leq H^b$.

\section{Main Results}

In this section, we study the leader-follower flocking problem for networked Lagrange systems. The goal is to design $u_i$ for each follower to achieve the leader-follower flocking. That is, the followers move cohesively with the leader (connectivity maintenance) and avoid collisions for all time and eventually achieve velocity matching with the leader ($||\dot{q}_i(t)-\dot{q}_0(t)|| \to 0$) in the presence of unknown parameters under only local interaction defined by the proximity graph. Before moving on, the following auxiliary variables are defined:
\begin{equation}\label{equ2}
s_i=\dot{q}_i-v_i,\qquad
\tilde{q}_i=q_i-q_0,\qquad
\tilde{v}_i=v_i-\dot{q}_0,
\end{equation}
where $v_i$ is agent $i$'s estimate of the leader's velocity to be designed later. Note that
\begin{equation}\label{sl-var}
s_i=\dot{\tilde{q}}_i-\tilde{v}_i.
\end{equation}

\subsection{Flocking when the leader has a constant velocity}
In this subsection, we consider the case where the leader has a constant velocity. We propose the following distributed control algorithm
\begin{align}
u_i=&\hat{u}_i+Y_i(q_i,\dot{q}_i,\dot{v}_i,v_i)\hat{\theta}_i,\label{input-cons-vel}
\\
\hat{u}_i=&-\sum_{j=0}^n \frac{\partial V_{ij}}{\partial q_i}-\gamma\sum_{j=0}^{n} a_{ij}(t) (\dot{q}_i-\dot{q}_j ),\label{u-hat-cons-vel}
\\
\dot{v}_i=&-\sum_{j=0}^n \frac{\partial V_{ij}}{\partial q_i}-\gamma\sum_{j=0}^{n} a_{ij}(t) (\dot{q}_i-\dot{q}_j ),\label{v-hat-cons-vel}
\\
\dot{\hat{\theta}}_i=&-\Gamma_i Y_i^T(q_i,\dot{q}_i,\dot{v}_i,v_i) s_i,\label{par-estm-hat-cons-vel}
\end{align}
where $a_{ij}(t)$ is the edge weight associated with the proximity graph $\overline{G}$ defined in Section II-B, $V_{ij}$ is the potential function between agents $i$ and $j$ to be designed, $\hat{\theta}_i$ is the estimate of the unknown but constant parameter $\theta_i$, $s_i$ is defined in \eqref{equ2}, $\gamma$ is a positive constant, and $\Gamma_i$ is a symmetric positive-definite matrix representing the adaptation gain.

\begin{remark}
Here $v_i$ is the reference velocity, which introduces the partial derivatives of the potential functions in the estimators and it is a key to our problem. It is worthy mentioning that \eqref{v-hat-cons-vel} has a similar form of the reference velocity derivative proposed in \cite{Wang13}, where the partial derivatives are replaced by the position synchronization term. Compared to the position and velocity synchronization problem consider in \cite{Wang13}, here we study the flocking problem (connectivity maintenance, collision avoidance, and velocity matching with a moving leader whose information is available to only the followers in its proximity).
\end{remark}

The potential function $V_{ij}$ is defined as follows (see \cite{cao2012distributed})
\begin{enumerate}
\item When $||q_i(0)-q_j(0)|| \ge R$, $V_{ij}$ is a differentiable nonnegative function of $||q_i-q_j||$ satisfying the conditions:
\begin{itemize}
\item[i)] $V_{ij}=V_{ji}$ achieves its unique minimum when $||q_i-q_j ||$ is equal to the value $\overline{d}_{ij}$, where $\overline{d}_{ij}<R$.
\item[ii)] $V_{ij} \to \infty$ as $||q_i-q_j ||\to 0$.
\item[iii)] $\frac{\partial V_{ij}}{\partial (||q_i-q_j||)}=0$ if $||q_i-q_j ||\ge R$.
\item[iv)] $V_{ii}=c, i=1,\ldots,n$, where $c$ is a positive constant.
\end{itemize}
\item When $||q_i(0)-q_j(0)||<R$, $V_{ij}$ is defined as above except that condition iii) is replaced with the condition that $V_{ij} \to \infty$ as $||q_i-q_j|| \to R$.
\end{enumerate}
The motivation of $V_{ij}$ is to maintain the initial connectivity pattern and to avoid collision.

In the control algorithm \eqref{input-cons-vel}-\eqref{par-estm-hat-cons-vel}, the term $-\sum_{j=0}^{n} \frac{\partial V_{ij}}{\partial q_i}$ is used for collision avoidance and connectivity maintenance while the term $- \sum_{j=0}^{n} a_{ij}(t)(\dot{q}_i-\dot{q}_j)$ is used for velocity matching. The control algorithm \eqref{input-cons-vel}-\eqref{par-estm-hat-cons-vel} is distributed in the sense that each agent uses only its own position and velocity and the relative position and relative velocity between itself and its neighbors.

\begin{theorem}\label{thm:ELflocking-Cons}
Suppose that at the initial time $t=0$, the leader has directed paths to all followers and there is no collision among the agents. Using \eqref{input-cons-vel}-\eqref{par-estm-hat-cons-vel} for \eqref{EL-dyn}, the leader-follower flocking is achieved.
\end{theorem}

\emph{Proof}: By using the property (P3) of the Lagrange dynamics \eqref{EL-dyn}, it follows that
${M}_i(q_i) \dot{v}_i+{C}_i(q_i,\dot{q}_i) v_i+{g}_i(q_i)=Y_i (q_i,\dot{q}_i,\dot{v}_i,v_i){\theta}_i$.
Then using \eqref{EL-dyn}, \eqref{equ2} and \eqref{input-cons-vel}, we have the following closed-loop system
\begin{align}\label{closed-loop-cons-vel}
M_i (q_i ) \dot{s}_i+C_i (q_i,\dot{q}_i ) s_i&=\hat{u}_i-Y_i(q_i,\dot{q}_i,\dot{v}_i,v_i)\tilde{\theta}_i,
\end{align}
where $\tilde{\theta}_i=\theta_i-\hat{\theta}_i$.
We first define the following non-negative function, which is a common Lyapunov function candidate used in the literature \cite{mei2012distributed,NunoOBH11_TAC,Wang13,Wang14_TAC} with different definition of $s_i$
\begin{equation}\label{V-1-fnc}
V_1=\frac{1}{2} \sum_{i=1}^n s_i^T M_i(q_i) s_i+\frac{1}{2} \sum_{i=1}^n \tilde{\theta}_i^T \Gamma_i^{-1} \tilde{\theta}_i.
\end{equation}
The derivative of $V_1$ is given as
\begin{align}\label{V1-dot-cons-vel}
\dot{V}_1 =&\sum_{i=1}^n [s_i^T M_i(q_i) \dot{s}_i+\frac{1}{2} s_i^T \dot{M}_i(q_i) s_i-\tilde{\theta}_i^T \Gamma_i^{-1} \dot{\hat{\theta}}_i] \notag \\
=&\sum_{i=1}^n s_i^T \hat{u}_i,
\end{align}
where we have used 
the property (P2) and \eqref{par-estm-hat-cons-vel} to obtain the last equality.
To maintain the
initial connectivity pattern and to avoid collision, we then define the following negative function by the combination of the potential functions
\allowdisplaybreaks[4]
\begin{align}\label{V2-fnc}
V_2=\frac{1}{2} \sum_{i=1}^n \sum_{j=1}^n V_{ij}+\sum_{i=1}^n V_{i0}.
\end{align}
Its derivative can be written as
\allowdisplaybreaks[4]
\begin{align}
\dot{V}_2=&\frac{1}{2} \sum_{i=1}^n \sum_{j=1}^n (\dot{q}_i^T \frac{\partial V_{ij}}{\partial q_i }+\dot{q}_j^T \frac{\partial V_{ij}}{\partial q_j})+\sum_{i=1}^n (\dot{q}_i^T \frac{\partial V_{i0}}{\partial q_i }+\dot{q}_0^T \frac{\partial V_{i0}}{\partial q_0 })\notag\\
=&\sum_{i=1}^n \sum_{j=1}^n \dot{q}_i^T \frac{\partial V_{ij}}{\partial q_i }+\sum_{i=1}^n (\dot{q}_i^T \frac{\partial V_{i0}}{\partial q_i }-\dot{q}_0^T \frac{\partial V_{i0}}{\partial q_i })\notag\\
=&\sum_{i=1}^n \sum_{j=1}^n (\dot{q}_i-\dot{q}_0)^T \frac{\partial V_{ij}}{\partial q_i }+\sum_{i=1}^n \dot{\tilde{q}}_i^T \frac{\partial V_{i0}}{\partial q_i } \notag \\
=&\sum_{i=1}^n \sum_{j=0}^n \dot{\tilde{q}}_i^T \frac{\partial V_{ij}}{\partial q_i},\notag
\end{align}
where we have used Lemma 3.1 in \cite{cao2012distributed} and the fact that $\frac{\partial V_{ij}}{\partial q_i }=-\frac{\partial V_{ij}}{\partial q_j}$ to obtain the second equality, and have used the fact that $\sum_{i=1}^n \sum_{j=1}^n \dot{q}_0^T \frac{\partial V_{ij}}{\partial q_i}=\dot{q}_0^T\sum_{i=1}^n \sum_{j=1}^n  \frac{\partial V_{ij}}{\partial q_i}=0$ to obtain the third equality.

Now consider the following Lyapunov function candidate
\begin{equation}\label{lyp-fnc-cons-vel}
V=V_1+\frac{1}{2} \sum_{i=1}^n \tilde{v}_i^T \tilde{v}_i +V_2.
\end{equation}
Then the derivative of $V$ is given as
\begin{align}
\dot{V}=&\sum_{i=1}^n s_i^T \hat{u}_i+\sum_{i=1}^n \tilde{v}_i^T \dot{\tilde{v}}_i+\sum_{i=1}^n \sum_{j=0}^n \dot{\tilde{q}}_i^T \frac{\partial V_{ij}}{\partial q_i}.\notag
\end{align}
Since the leader's velocity $\dot{q}_0$ is constant, we have $\dot{\tilde{v}}_i=\dot{v}_i=\hat{u}_i$ according to \eqref{u-hat-cons-vel} and \eqref{v-hat-cons-vel}. It follows that
\begin{align}\label{V-dot-cons-vel}
\dot{V}=&\sum_{i=1}^n \dot{\tilde{q}}_i^T \hat{u}_i+\sum_{i=1}^n \sum_{j=0}^n \dot{\tilde{q}}_i^T \frac{\partial V_{ij}}{\partial \tilde{q}_i} \notag \\
=&-\sum_{i=1}^n \sum_{j=0}^n \gamma a_{ij}(t) \dot{\tilde{q}}_i^T (\dot{\tilde{q}}_i-\dot{\tilde{q}}_j ),
\end{align}
where we have used \eqref{sl-var} to obtain the first equality and have used \eqref{u-hat-cons-vel} and $\dot q_i-\dot q_j=\dot{\tilde{q}}_i-\dot{\tilde{q}}_j$ to obtain the second equality.
Eq. \eqref{V-dot-cons-vel} can be written in a compact form as
\allowdisplaybreaks[4]
\begin{equation}\label{lyp-fnc-derivative-const-vel}
\dot{V}=-\gamma\dot{\tilde{q}}^T [H(t) \otimes I_p]\dot{\tilde{q}},
\end{equation}
where $\tilde{q}$ is a column stack vector of $\tilde{q}_i$, $i=1,\ldots,n$, and $H(t)$ is the leader-follower topology matrix at time $t$ defined in Section II-B. Note that $H(t)$ is symmetric positive semi-definite. It follows that $\dot{V}$ is negative semi-definite. Therefore, from $V \ge 0$ and $\dot{V} \le 0$, it can be concluded that $V$ is bounded and thus $s_i$, $\tilde{\theta}_i$, $\tilde{v}_i$, $V_{ij}\in\mathbb{L}_{\infty}$. Since $V_{ij}$ is bounded, it is guaranteed that there is no collision and no edge in the graph $\overline{G}(0)$ will be lost.
In other words, for any pair of agents $i$, $j$, there exist positive constants $0<R_{\min}\leq R_{\max}<R$, such that
\begin{align}\label{equ101}
&\|q_i(t)-q_j(t)\|\in[R_{\min}, R_{\max}], \quad\quad \text{if} ~\|q_i(0)-q_j(0)\|<R,\notag\\
&\|q_i(t)-q_j(t)\|\in[R_{\min}, (n-1)R_{\max}],\quad \text{otherwise.}
\end{align}
Hence, we can conclude that the graph $\overline{G}(0)$ is a subgraph of the graph $\overline{G}(t)$ for all $t\geq 0$. It follows from Lemma \ref{lemma:psd} that $H(0) \leq H(t)$. Therefore, we can get from \eqref{lyp-fnc-derivative-const-vel} that
\begin{align}\label{equ100}
\dot{V} \leq -\gamma \dot{\tilde{q}}^T[H(0) \otimes I_p] \dot{\tilde{q}}.
\end{align}
Since in $\overline{G}(0)$ the leader has directed paths to all followers, it follows from Lemma \ref{lemma:tree} that $H(0)$ is symmetric positive definite.
Integrating both sides of \eqref{equ100}, we can obtain that $\dot{\tilde q}\in\mathbb{L}_2$.
Note that $\dot{q}_0$ is constant and hence bounded. Combining the above boundedness arguments we can get from \eqref{equ2} that $\dot{q}_i$, $\dot{\tilde q}_i$, $v_i\in\mathbb{L}_{\infty}$.
Since $V_{ij}$ is continuously differentiable, we can get from \eqref{equ101} that $\frac{\partial V_{ij}}{\partial q_i}\in\mathbb{L}_{\infty}$. From \eqref{u-hat-cons-vel} and \eqref{v-hat-cons-vel}, we have
$\hat u_i, \dot v_i\in\mathbb{L}_{\infty}$.
Then from \eqref{closed-loop-cons-vel} and the property (P1), it can be concluded that $\dot{s}_i\in\mathbb{L}_{\infty}$. By noting that $\dot{s}_i=\ddot{q}_i-\dot{v}_i$, it follows that $\ddot{q}_i\in\mathbb{L}_{\infty}$. Overall, we have $\dot{\tilde{q}}_i\in\mathbb{L}_{\infty}\bigcap\mathbb{L}_2$ and $\ddot{\tilde{q}}_i\in\mathbb{L}_{\infty}$. From Barbalat's lemma \cite{slotine1991applied}, we can conclude that $\dot{\tilde{q}}_i \to 0$, that is, $||\dot{q}_i-\dot{q}_0|| \to 0$ asymptotically.
\endproof

\begin{remark} As it can be seen, by using the control law \eqref{input-cons-vel}-\eqref{par-estm-hat-cons-vel} for \eqref{EL-dyn}, the followers can track the leader with the same velocity while avoiding collision and maintaining the initial connectivity. Note that with our algorithm design, as long as at the initial time the connectivity is maintained and there is no collision, the connectivity maintenance and collision avoidance are ensured for all time. The proposed algorithm is continuous and accounts for unknown parameters of the agents' dynamics. 
\end{remark}

\subsection{Flocking when the leader has a varying velocity}
In this subsection, we consider the case when the leader moves with a varying velocity. In this case, the problem is more difficult to tackle since all followers must track the leader while the leader's velocity changes over time and the leader is a neighbor of only a subset of the followers in its proximity. In the remainder of the paper, we have the following assumption on the leader.
\begin{assumption}\label{assume:leader}
The leader's velocity $\dot q_0$ and acceleration $\ddot q_0$ are both bounded. It is assumed that $||\mathbf{1}_n \otimes \ddot{q}_0|| \le \sigma_l$, where $\sigma_l$ is a positive constant.
\end{assumption}

We propose the following distributed control algorithm
\begin{align}
u_i=&\hat{u}_i+Y_i(q_i,\dot{q}_i,\dot{v}_i,v_i)\hat{\theta}_i,\label{input-varying-vel}
\\
\hat{u}_i=&-\sum_{j=0}^n \frac{\partial V_{ij}}{\partial q_i } -\alpha \sum_{j=0}^n a_{ij}(t) \mbox{sgn} (\dot{q}_i-\dot{q}_j) -\alpha \mbox{sgn}(s_i),\label{u-hat-varying-vel}
\\
\dot{v}_i=&-\sum_{j=0}^n \frac{\partial V_{ij}}{\partial q_i }- \alpha \sum_{j=0}^n a_{ij}(t)  \mbox{sgn}(\dot{q}_i-\dot{q}_j), \label{v-dot-varying-vel}
\\
\dot{\hat{\theta}}_i=&-\Gamma_i Y_i^T(q_i,\dot{q}_i,\dot{v}_i,v_i) s_i, \label{par-estm-varying-vel}
\end{align}
where $\alpha$ is a positive constant, and $a_{ij}(t)$, $V_{ij}$, $s_i$ and $\Gamma_i$ are defined as in Section III-A.

\begin{theorem}\label{thm:ELflocking-Vary}
Suppose that at the initial time $t=0$, the leader has directed paths to all followers and there is no collision among the agents. Using  \eqref{input-varying-vel}-\eqref{par-estm-varying-vel} for \eqref{EL-dyn}, if $\alpha > \max \{ \sigma_l, \frac{\sigma_l }{\sqrt{\lambda_{\min} [H(0)]}} \}$\footnote{Since at the initial time $t=0$, the leader has directed paths to all followers, we can get from Lemma \ref{lemma:tree} that $\lambda_{\min}[H(0)]>0$, and thus the term $\frac{\sigma_l }{\sqrt{\lambda_{\min} [H(0)]}}$ is well defined.}, then the leader-follower flocking is achieved.
\end{theorem}

\emph{Proof}: Consider the same Lyapunov function candidate $V$ defined in \eqref{lyp-fnc-cons-vel}. Note that using \eqref{input-varying-vel} for \eqref{EL-dyn}, where $\hat{u}_i$ is given by \eqref{u-hat-varying-vel}, both \eqref{closed-loop-cons-vel} and \eqref{V1-dot-cons-vel} still hold. The derivative of $V$ is given as
\begin{align}\notag
\dot{V}=\sum_{i=1}^n s_i^T \hat{u}_i+\sum_{i=1}^n \tilde{v}_i^T \dot{\tilde{v}}_i+\sum_{i=1}^n \sum_{j=0}^n \dot{\tilde{q}}_i^T \frac{\partial V_{ij}}{\partial {q}_i}.
\end{align}
Note from $\eqref{u-hat-varying-vel}$ and $\eqref{v-dot-varying-vel}$ that
\begin{align}\notag
\dot{v}_i=\hat{u}_i+\alpha \mbox{sgn}(s_i).
\end{align}
Also note from $\eqref{equ2}$ that $\dot{\tilde{v}}_i=\dot{v}_i-\ddot{q}_0$.
It follows that
\begin{align}\notag
\dot{V}=&\sum_{i=1}^n s_i^T \hat{u}_i+\sum_{i=1}^n \tilde{v}_i^T[\hat{u}_i+\alpha \mbox{sgn}(s_i)-\ddot{q}_0] \notag \\
&+\sum_{i=1}^n\sum_{j=0}^n \dot{\tilde{q}}_i^T \frac{\partial V_{ij}}{\partial {q}_i }.
\end{align}
Note from \eqref{sl-var} that $\tilde{v}_i=\dot{\tilde{q}}_i-s_i$. Therefore, it follows that
\begin{align}\label{equ50}
\dot{V}=&\sum_{i=1}^n \dot{\tilde{q}}_i^T[\hat{u}_i+\alpha \mbox{sgn}(s_i)-\ddot{q}_0]-\sum_{i=1}^n s_i^T[\alpha \mbox{sgn}(s_i)-\ddot{q}_0] \notag \\
&+\sum_{i=1}^n \sum_{j=0}^n \dot{\tilde{q}}_i^T \frac{\partial V_{ij}}{\partial {q}_i }.
\end{align}
Substituting $\hat{u}_i$ defined in \eqref{u-hat-varying-vel} to \eqref{equ50}, we can get	
\begin{align*}
\dot{V}=&-\alpha  \dot{\tilde{q}}^T [D_F(t) \otimes I_p] \mbox{sgn}
([D^T_F(t) \otimes I_p]\dot{\tilde{q}}) \\
& - \alpha \dot{\tilde{q}}^T[\Lambda (t) \otimes I_p] \mbox{sgn}
([\Lambda (t) \otimes I_p]\dot{\tilde{q}}) -\dot{\tilde{q}}^T (\mathbf{1}_n \otimes \ddot{q}_0) \notag \\
&-\alpha ||s||_1+s^T (\mathbf{1}_n \otimes \ddot{q}_0) \\
=&-\alpha ||[D^T_F(t) \otimes I_p]\dot{\tilde{q}}||_1 -\alpha ||[\Lambda (t) \otimes I_p]\dot{\tilde{q}}||_1 \\
& -\dot{\tilde{q}}^T (\mathbf{1}_n \otimes \ddot{q}_0)-\alpha ||s||_1 +s^T (\mathbf{1}_n \otimes \ddot{q}_0) \\
\le& -\alpha ||[D^T_F(t) \otimes I_p]\dot{\tilde{q}}|| -\alpha ||[\Lambda (t) \otimes I_p]\dot{\tilde{q}}|| \\
&+||\mathbf{1}_n \otimes \ddot{q}_0||\cdot||\dot{\tilde{q}}||-\alpha ||s||+||\mathbf{1}_n \otimes \ddot{q}_0||\cdot ||s||
\end{align*}
where $s$ and $\tilde{q}$ are, respectively, the column stack vectors of all $s_i$'s and $\tilde{q}_i$'s, $i=1,\ldots,n$, and we have used the fact that $\|\cdot\|\leq \|\cdot\|_1$ for any vector to obtain the inequality. Since $||\mathbf{1}_n \otimes \ddot{q}_0|| \le \sigma_l$, we have
\begin{align*}
\dot{V} \le& -\alpha ||
	\begin{bmatrix}
	[D^T_F(t) \otimes I_p] \\
	[\Lambda (t) \otimes I_p]
	\end{bmatrix}	\dot{\tilde{q}}|| +\sigma_l||\dot{\tilde{q}}|| -\alpha ||s||+\sigma_l ||s|| \\
=& -\alpha \sqrt{\dot{\tilde{q}}^T
		\begin{bmatrix}
		[D_F(t) \otimes I_p] && [\Lambda (t) \otimes I_p]
		\end{bmatrix}
		\begin{bmatrix}
		[D^T_F(t) \otimes I_p] \\
		[\Lambda (t) \otimes I_p]
		\end{bmatrix} \dot{\tilde{q}}} \\
&+\sigma_l||\dot{\tilde{q}}||-(\alpha -\sigma_l) ||s|| \\
=& -\alpha \sqrt{\dot{\tilde{q}}^T
		[H(t) \otimes I_p] \dot{\tilde{q}}} +\sigma_l||\dot{\tilde{q}}|| -(\alpha -\sigma_l) ||s|| \\
\leq & -\alpha \sqrt{\lambda_{\min}[H(t)]}||\dot{\tilde{q}}|| +\sigma_l||\dot{\tilde{q}}|| -(\alpha -\sigma_l) ||s||
\end{align*}
where we have used the equation $D_F(t)D^T_F(t)+\Lambda^2 (t)=L_F(t)+\Lambda(t)=H(t)$ to obtain the second equality and have used the fact that $H(t)$ is positive semi-definite to obtain the last inequality. 
Note that at the initial time $t=0$, the leader has directed paths to all followers. We can get from Lemma \ref{lemma:tree} that $H(0)$ is symmetric positive definite and thus $\lambda_{\min}[H(0)]>0$.
Since $\alpha > \max \{ \sigma_1, \frac{\sigma_1}{\sqrt{\lambda_{\min} [H(0)]}} \}$, we have at time $t=0$ that,
 \begin{align}\label{equ25}
	\dot{V}(t) \le -(\alpha \sqrt{\lambda_{\min}[H(0)]} - \sigma_l)||\dot{\tilde{q}}||  \leq 0.
	\end{align}
Note that although the control input $u_i$ is discontinuous, the positions of the agents are continuous and $H(t)$ changes according to the relative positions among the agents. If $H(t)$ changes at some time, there exists $t_1>0$ such that, ${H(t)=H(0)}$ for $t\in[0,t_1)$ and $H(t_1)\neq H(0)$. Therefore, we have
\begin{equation}\notag
\dot{V}(t) \le -(\alpha \sqrt{\lambda_{\min}[H(0)]} - \sigma_l)||\dot{\tilde{q}}|| \leq 0, \qquad t\in[0,t_1),
\end{equation}
which implies that for $t\in[0,t_1)$, $V_{ij}\in\mathbb{L}_{\infty}$ for all pairs of $q_i(t)$ and $q_j(t)$. Since $V_{ij}$ is continuous, we can conclude that $V_{ij}\in\mathbb{L}_{\infty}$ when $t=t_1$. From the definition of $V_{ij}$, it follows that there is no collision and also no edge in the
graph $\overline G(0)$ will be lost for $t\in[0,t_1]$. Therefore, the only possibility that $H(t)$ changes at $t=t_1$ is that, some edges are added in the graph. It implies that $\overline G(0)$ is a subgraph of $\overline G(t_1)$. We can then get from Lemma \ref{lemma:psd} that $H(0)\leq H(t_1)$ and thus $\lambda_{\min}[H(0)]\leq \lambda_{\min}[H(t_1)]$. Therefore, at time $t=t_1$,
\begin{align*}
\dot{V}(t) \le& -(\alpha \sqrt{\lambda_{\min} [H(t_1)]} -\sigma_l)||\dot{\tilde{q}}|| \\
\le& -(\alpha \sqrt{\lambda_{\min} [H(0)]}-\sigma_l)||\dot{\tilde{q}}||.
\end{align*}
Following the same argument, if $H(t)$ changes at ${t=t_i>t_1}$, $i=2,\ldots$, we can get that $V_{ij}$ will always be bounded. Hence there is no collision and no edge in the
graph $\overline G(0)$ will be lost.
This in turn implies that for all $t \in [t_i,t_{i+1})$, $\overline{G}(0)$ is a subgraph of $\overline{G}_i(t)$.
It thus follows that for all $t \in [t_i,t_{i+1})$, $H(0) \leq H_i(t)$ and $\lambda_{\min} [H(0)] \leq \lambda_{\min} [H_i(t)]$.
That is, $\alpha > \frac{\sigma_1}{\sqrt{\lambda_{\min} [H(0)]}} \geq \frac{\sigma_1}{\sqrt{\lambda_{\min} [H_i(t)]}}$ for all $t \geq 0$. Hence \eqref{equ25} holds for all $t \geq 0$.
We then can get that $s_i$, $\tilde{\theta}_i$, $\tilde{v}_i\in\mathbb{L}_{\infty}$.
Since $V(t) \geq 0$ and $\dot{V}(t) \leq 0$, it is concluded that $V_{\infty} \defeq \lim_{t \to \infty} V(t)\in [0,V(0)]$ exists. Thus, integrating both sides of \eqref{equ25}, we can obtain that
$\dot{\tilde q}\in \mathbb{L}_1$. Note from Assumption \ref{assume:leader} that both $\dot{q}_0$ and $\ddot q_0$ are bounded. Combining the above boundedness arguments, we can get from \eqref{equ2} that $\dot{q}_i$, $\dot{\tilde q}_i$, $v_i\in\mathbb{L}_{\infty}$.
Following the same statements from the proof of Theorem \ref{thm:ELflocking-Cons}, we can conclude that $\frac{\partial V_{ij}}{\partial q_i}\in\mathbb{L}_{\infty}$. From \eqref{u-hat-varying-vel} and \eqref{v-dot-varying-vel}, we have
$\hat u_i$, $\dot v_i\in\mathbb{L}_{\infty}$.
Then from the closed-loop dynamics for each follower and (P1), we have $\dot{s}_i\in\mathbb{L}_{\infty}$. By noting that $\dot{s}_i=\ddot{q}_i-\dot{v}_i$, it follows that $\ddot{q}_i\in\mathbb{L}_{\infty}$ and thus $\ddot{\tilde q}_i\in\mathbb{L}_{\infty}$. Overall, we have $\dot{\tilde{q}}_i\in\mathbb{L}_{\infty}\bigcap\mathbb{L}_1$ and $\ddot{\tilde{q}}_i\in\mathbb{L}_{\infty}$. From Barbalat's lemma \cite{slotine1991applied}, we can conclude that $\dot{\tilde{q}}_i \to 0$. That is, $||\dot{q}_i-\dot{q}_0|| \to 0$ asymptotically. \endproof


\begin{remark} As it can be seen, the proposed algorithm \eqref{input-varying-vel}-\eqref{par-estm-varying-vel} guarantees that the leader-follower flocking is achieved when the leader has a varying velocity in the presence of unknown parameters. Therefore, despite the hard restrictions such as nonlinear Lagrange dynamics, unknown models' parameters, and the existence of a moving leader with a varying velocity, the control input \eqref{input-varying-vel}-\eqref{par-estm-varying-vel} solves the flocking problem.
\end{remark}

\begin{remark}
Due to the existence of the signum function, the closed-loop dynamics of \eqref{EL-dyn} using \eqref{input-varying-vel} is discontinuous. The solution should be investigated in terms of differential inclusions. Note that the signum function is measurable and locally essentially bounded. Therefore, from the nonsmooth analysis in \cite{filippov1960differential}, the Filippov solutions for the closed-loop dynamics always exist. Because the Lyapunov function candidate in the proof of Theorem \ref{thm:ELflocking-Vary} is continuously differentiable and the set-valued Lie derivative of the Lyapunov function is a singleton at the discontinuous point, the proof of Theorem \ref{thm:ELflocking-Vary} still holds. To avoid symbol redundancy, we do not use the differential inclusions in the proof. It is worthy mentioning that the drawback of the signum function is the potential chattering behavior. In practice, a simple and useful way to avoid the discontinuous of the control action is to replace the signum function by a smooth function such as $\tanh(\cdot)$, with which satisfactory performance can still be achieved, as confirmed in our later simulation.
\end{remark}

\begin{remark}
The case of a leader with a constant velocity is a special case of a leader with a varying velocity. Hence we can also use the algorithm \eqref{input-varying-vel}-\eqref{par-estm-varying-vel} for the leader-follower flocking problem when the leader has a constant velocity. However, the algorithm \eqref{input-cons-vel}-\eqref{par-estm-hat-cons-vel} is continuous. 
In contrast, the algorithm \eqref{input-varying-vel}-\eqref{par-estm-varying-vel} is discontinuous and may cause the chattering issues. 
Therefore, when the leader has a constant velocity, the algorithm \eqref{input-cons-vel}-\eqref{par-estm-hat-cons-vel} is more favorable than the algorithm \eqref{input-varying-vel}-\eqref{par-estm-varying-vel}.
\end{remark}

\subsection{Fully distributed flocking when the leader has a varying velocity}
In the previous section, all agents use common gains in their control inputs and the gains should be above certain bounds which are actually determined by the global information $(\lambda_{\min} [H(0)],\sigma_l)$. Therefore, the algorithm \eqref{input-varying-vel}-\eqref{par-estm-varying-vel} is not fully distributed. In this section, the previous algorithm is extended to be fully distributed and gain adaptation laws are introduced. Here the control algorithm for each follower is designed as
\begin{align}
	u_i&=\hat{u}_i+Y_i(q_i,\dot{q}_i,\dot{v}_i,v_i)\hat{\theta}_i,\label{equ28}
	\\
 \hat{u}_i &=-\sum_{j=0}^n \frac{\partial V_{ij}}{\partial q_i } -\sum_{j=0}^n \alpha_{ij} a_{ij}(t)  \mbox{sgn} (\dot{q}_i-\dot{q}_j)-\beta_{i} \mbox{sgn}(s_i),\label{equ29}
	\\
	\dot{v}_i &=-\sum_{j=0}^n \frac{\partial V_{ij}}{\partial q_i }-\sum_{j=0}^n \alpha_{ij} a_{ij}(t)  \mbox{sgn} (\dot{q}_i-\dot{q}_j),\label{equ30}
	\\
	\dot{\alpha}_{ij} &= \gamma_{1i}a_{ij}(t) ||\dot{{q}}_i-\dot q_0||_1,
	\label{equ32}
	\\
	\dot{\beta}_i &= \gamma_{2i}||s_i||_1 , \label{equ33}
	\\	\dot{\hat{\theta}}_i &=-\Gamma_i Y_i^T(q_i,\dot{q}_i,\dot{v}_i,v_i) s_i,\label{equ34}
	\end{align}
where $a_{ij}(t)$, $V_{ij}$, $s_i$ and $\Gamma_i$ are defined in Section III-A, $\gamma_{1i}$, $\gamma_{2i}$ are positive constants, and $\alpha_{ij}(t)$, $\beta_{i}(t)$ are varying gains with $\alpha_{ij}(0), \beta_{i}(0)\geq 0$.

\begin{remark}
The gain adaptation laws \eqref{equ32} and \eqref{equ33} are inspired by recent results on adaptive gain design for multi-agent systems \cite{YuRenZCL13_automatica,LiRenLiuXie13_automatica,MeiRenChenMa13_automatica}.
 The intuition behind \eqref{equ32} and \eqref{equ33} is that the control gains in Theorem \ref{thm:ELflocking-Vary} must be above certain lower bounds. Under \eqref{equ32} and \eqref{equ33}, as long as the velocity matching is not achieved, the gains will always increase, eventually rendering the agents to achieve velocity matching. The drawback of \eqref{equ32} and \eqref{equ33} is that the 1-norm of the signals will result in the non-stop increase of the gains in the presence of disturbances or measurement errors.
Here we just show the theoretical analysis in the ideal situation. In practice, one alteration is to introduce a small bound on the right hand sides (RHSs) of \eqref{equ32} and \eqref{equ33}. When the RHSs of \eqref{equ32} and \eqref{equ33} are within some given bound, $\alpha_{ij}$ and $\beta_i$ stops increasing. 
\end{remark}

\begin{theorem}\label{thm:distributed-Flock}
Suppose that at the initial time $t=0$, the leader has directed paths to all followers and there is no collision among the agents. Using \eqref{equ28}-\eqref{equ34} for \eqref{EL-dyn}, the leader-follower flocking is achieved.
\end{theorem}

\emph{Proof}: Define $V_3=\sum_{i=1}^n\frac{1}{4\gamma_{1i}} \sum_{j=1}^n (\alpha_{ij}-\bar{\alpha})^2+\sum_{i=1}^n \frac{1}{2\gamma_{1i}}(\alpha_{i0}-\bar{\alpha})^2+\sum_{i=1}^n\frac{1}{2\gamma_{2i}} (\beta_i-\bar{\beta})^2$, where $\bar{\alpha}$ and $\bar{\beta}$ are chosen such that $\bar{\alpha} > \frac{\sigma_l }{\sqrt{\lambda_{\min}[H(0)]}}$ and $\bar{\beta}>\sigma_l$. The derivative of $V_3$ is given as
 \begin{align}\label{equ37}
	\dot{V}_3
	=&\sum_{i=1}^n \sum_{j=1}^n \frac{1}{2} a_{ij}(t)(\alpha_{ij}-\bar{\alpha}) ||\dot{q}_i-\dot{q}_j||_1 \notag \\
	&+\sum_{i=1}^n a_{i0}(t)(\alpha_{i0}-\bar{\alpha})||\dot{\tilde{q}}_i||_1 +\sum_{i=1}^n (\beta_i-\bar{\beta})||s_i||_1 \notag \\
	=&\sum_{i=1}^n \sum_{j=1}^n \frac{1}{2} a_{ij}(t)(\alpha_{ij}-\bar{\alpha}) (\dot{q}_i-\dot{q}_j)^T \mbox{sgn}(\dot{q}_i-\dot{q}_j)\notag \\
	&+\sum_{i=1}^na_{i0}(t)(\alpha_{i0}-\bar{\alpha})(\dot{q}_i-\dot{q}_0)^T \mbox{sgn}(\dot{q}_i-\dot{q}_0) \notag \\
	&+\sum_{i=1}^n (\beta_i-\bar{\beta})||s_i||_1 \notag \\
	=&\sum_{i=1}^n \sum_{j=1}^n  a_{ij}(t)(\alpha_{ij}-\bar{\alpha}) \dot{\tilde{q}}^T_i \mbox{sgn}(\dot{q}_i-\dot{q}_j)\notag \\
	&+\!\sum_{i=1}^na_{i0}(t)(\alpha_{i0}\!-\!\bar{\alpha})\dot{\tilde{q}}^T_i \mbox{sgn}(\dot{q}_i\!-\!\dot{q}_0)\!+\!\sum_{i=1}^n (\beta_i\!-\!\bar{\beta})||s_i||_1 \notag \\
	=&\sum_{i=1}^n [\sum_{j=0}^n  a_{ij}(t)(\alpha_{ij}\!-\!\bar{\alpha}) \dot{\tilde{q}}^T_i \mbox{sgn}(\dot{\tilde{q}}_i\!-\!\dot{\tilde{q}}_j)\!+\!(\beta_i\!-\!\bar{\beta})||s_i||_1].
	\end{align}
Now we introduce the following Lyapunov function candidate
\begin{align}
V=&V_1+\frac{1}{2} \sum_{i=1}^n \tilde{v}_i^T \tilde{v}_i+V_2+V_3,\notag
\end{align}
where $V_1$ is defined in \eqref{V-1-fnc} and $V_2$ is defined in \eqref{V2-fnc}. Following the proof of Theorem \ref{thm:ELflocking-Vary}, the derivative of $V$ is given as
\begin{align}\label{equ34}
\dot{V}=& \sum_{i=1}^n s_i^T \hat{u}_i+\sum_{i=1}^n \tilde{v}_i^T \dot{\tilde{v}}_i+\sum_{i=1}^n \sum_{j=0}^n \dot{\tilde{q}}_i^T \frac{\partial V_{ij}}{\partial q_i }+\dot{V}_3\notag\\
=& \sum_{i=1}^n s_i^T \hat{u}_i+\sum_{i=1}^n \tilde{v}_i^T [\hat{u}_i+\beta_{i} \mbox{sgn}(s_i)-\ddot{q}_0] \notag \\
&+\sum_{i=1}^n \sum_{j=0}^n \dot{\tilde{q}}_i^T \frac{\partial V_{ij}}{\partial q_i }+\dot{V}_3\notag\\
=& \sum_{i=1}^n \dot{\tilde{q}}_i^T [\hat{u}_i+\beta_{i} \mbox{sgn}(s_i)-\ddot{q}_0]-\sum_{i=1}^n s_i^T [\beta_{i} \mbox{sgn}(s_i)-\ddot{q}_0] \notag \\
&+\sum_{i=1}^n \sum_{j=0}^n \dot{\tilde{q}}_i^T \frac{\partial V_{ij}}{\partial q_i }+\dot{V}_3,
\end{align}
where we have used $\tilde{v}_i=\dot{\tilde{q}}_i-s_i$ to obtain the third equality.
Substituting \eqref{equ29} and \eqref{equ37} to \eqref{equ34} and doing some manipulation, we can get
\begin{align}
	\dot{V} \leq & -\bar{\alpha} ||[D^T_F(t) \otimes I_p]\dot{\tilde{q}}|| -\bar{\alpha} ||[\Lambda (t) \otimes I_p]\dot{\tilde{q}}|| \notag \\
	&+||\mathbf{1}_n \otimes \ddot{q}_0||\cdot||\dot{\tilde{q}}||-\bar{\beta} ||s|| +||\mathbf{1}_n \otimes \ddot{q}_0||\cdot ||s||,
	\end{align}
where $\tilde{q}$ and $s$ are, respectively, column stack vectors of $\tilde{q}_i$ and $s_i$, $i=1,\ldots,n$.
Since $||\mathbf{1}_n \otimes \ddot{q}_0 || \le \sigma_l$, we have
\begin{align*}
\dot{V} \leq & -\bar{\alpha} \sqrt{\lambda_{\min}[H(t)]} ||\dot{\tilde{q}}|| \notag  +\sigma_l ||\dot{\tilde{q}}||-\bar{\beta} ||s|| + \sigma_l ||s||.
\end{align*}
Again there exists $t_1 > 0$ such that, $H(t) = H(0)$ for $t \in [0,t_1)$. Since at time $t=0$ the parameters $\bar{\alpha}$ and $\bar{\beta}$ satisfy $\bar{\alpha} > \frac{\sigma_l }{\sqrt{\lambda_{\min} [H(0)]}}$ and $\bar{\beta} > \sigma_l$, we have
\begin{align*}
\dot{V}(t) \le -(\bar{\alpha} \sqrt{\lambda_{\min}[H(0)]} - \sigma_l)||\dot{\tilde{q}}|| \leq 0, \qquad t \in [0,t_1).
\end{align*}
Similar to the statements in the proof of Theorem \ref{thm:ELflocking-Vary}, it can be proved that for $t \in [0,t_1]$, $V_{ij}$ is bounded for all $i$, $j$ and there is no collision and also no edge in the graph $\overline{G}(0)$ will be lost. It can also be proved that $\dot{V}(t) \le 0$ for all $t \geq 0$. Then it can be concluded that $\dot{\tilde{q}}_i\in\mathbb{L}_{\infty}\bigcap\mathbb{L}_1$ and $\ddot{\tilde{q}}_i\in\mathbb{L}_{\infty}$. Hence from Barbalat's lemma, we conclude that $||\dot{q}_i-\dot{q}_0|| \to 0$ asymptotically.

\begin{remark}
In the fully distributed algorithm \eqref{equ28}-\eqref{equ32} adaptive gain schemes are introduced. This algorithm guarantees that the leader-follower flocking is achieved when the leader has a varying velocity in the presence of unknown parameters and there is no requirement of any global information. Therefore, despite the hard restrictions described in Section II-B and the existence of a moving leader with a varying velocity, the fully distributed control input \eqref{equ28}-\eqref{equ32} solves the flocking problem.
\end{remark}

\begin{remark}
In \cite{cao2012distributed}, the distributed flocking problem with a moving leader has been solved for multi-agent systems with single or double integrators. Here in this paper, we address the problem for networked nonlinear Lagrange systems with parametric uncertainties, which is more challenging. The algorithms in \cite{cao2012distributed} cannot deal with nonlinear Lagrange dynamics and account for fully distributed gain design. Besides, the algorithms in \cite{cao2012distributed} rely on both one-hop and two-hop neighbors' information, while only one-hop neighbors' information is required in our proposed algorithms. In \cite{mei2011distributed} a distributed coordinated tracking problem is studied for Lagrange systems with parametric uncertainties. However, the algorithms in \cite{mei2011distributed} cannot deal with the nonlinear flocking behavior or account for fully distributed gain design but still requires the two-hop neighbors' information. 
\end{remark}

\section{Simulation}

In this section, numerical simulation results are given to illustrate the effectiveness of the theoretical results obtained in Section III.
We consider the formation flying of four spacecraft, where the formation control is based on the relative translation with respect to a virtual point or chief spacecraft following a circular reference orbit \cite{AlfriendVHB_09}. The relative dynamics of the $i$th spacecraft is considered in a chief-fixed, LVLH rotating frame, which can be written as
\begin{align}
& m_i\ddot x_i-2 m_i n_0 \dot y_i- m_i n_0^2 x_i + \frac{m_i\mu_e(r_0+x_i)}{r_i^3}-\frac{m_i\mu_e}{r_0^2}= u_{ix},\notag\\
& m_i\ddot y_i+2 m_i n_0 \dot x_i-m_i n_0^2 y_i + \frac{m_i\mu_ey_i}{r_i^3} = u_{iy},\notag\\
& m_i\ddot z_i+ \frac{m_i\mu_ez_i}{r_i^3} = u_{iy},\notag
\end{align}
where $m_i$ is the unknown but constant mass of the $i$th spacecraft, $\mu_e$ is the gravitational constant of Earth, $r_0$ is the radius of the chief, $n_0=\sqrt{\mu_e/r_0^3}$ is the angular velocity of the reference orbit, $q_i\defeq [x_i, y_i, z_i]^T$ is the position of the $i$th spacecraft in the LVLH frame, and $u_i\defeq [u_{ix}, u_{iy}, u_{iz}]^T$ is the control input.
Let $M_i=m_iI_3$, $C_i=m_i\left(
                         \begin{array}{ccc}
                           0 & -2n_0 & 0 \\
                           2n_0 & 0 & 0 \\
                           0 & 0 & 0 \\
                         \end{array}
                       \right)
$, $g_i=m_i\left(
          \begin{array}{c}
            -n_0^2 x_i+\frac{\mu_e(r_0+x_i)}{r_i^3}-\frac{\mu_e}{r_0^2} \\
            -n_0^2 y_i+\frac{\mu_ey_i}{r_i^3}\\
            \frac{\mu_ez_i}{r_i^3} \\
          \end{array}
        \right)
$, then the relative translation dynamics can be written in the form of \eqref{EL-dyn}, with the unknown parameter $\theta_i=m_i$.

In the simulations, we let $m_i=30+5i~$kg, $r_0=7000~$km, and $R=200$~m. The initial positions of the leader and the four spacecraft are, respectively, $q_0(0)=[-80,200,0]^T$m, $q_1(0)=[-80, 90, 0]^T$m, $q_2(0)=[100, 90, 0]^T$m,
$q_3(0)=[100, -100, 0]^T$m, and $q_4(0)=[-80, -100, 0]^T$m. The initial velocities are assumed to be zero. The unique minimums of $V_{ij}$ are assumed to be $80$~m. Following \cite{cao2012distributed}, when $\|q_i(0)-q_j(0)\|\geq 200$, the potential functions are defined whose partial derivatives satisfy
\begin{align*}
& \pde{V_{ij}}{q_i} \\
& = \left\{
                    \begin{array}{ll}
                      0, & \|q_i-q_j\|>200; \\
                      \frac{(q_i-q_j)\cos(0.1\pi(\|q_i-q_j\|-80))}{250\|q_i-q_j\|} ,& 80<\|q_i-q_j\|\leq 200; \\
                      \frac{(q_i-q_j)(\|q_i-q_j\|-80)}{250\|q_i-q_j\|^2}, & \|q_i-q_j\|\leq 80.
                    \end{array}
                  \right.
\end{align*}
When $\|q_i(0)-q_j(0)\|< 200$, the potential functions are defined whose partial derivatives satisfy
\begin{align*}
 \pde{V_{ij}}{q_i} = \left\{
                    \begin{array}{ll}
                      \frac{(q_i-q_j)(\|q_i-q_j\|-80)}{25\|q_i-q_j\|(\|q_i-q_j\|-200)^2},& 80<\|q_i-q_j\|\leq 200; \\
                      \frac{(q_i-q_j)(\|q_i-q_j\|-80)}{250\|q_i-q_j\|^2}, & \|q_i-q_j\|\leq 80.
                    \end{array}
                  \right.
\end{align*}
%
In the first case, we simulate the case where the leader has a constant velocity under the control algorithm \eqref{input-cons-vel}-\eqref{par-estm-hat-cons-vel}. The constant velocity of the leader is assumed to be $\dot{q}_0=[0.1,0.1,0.2]^T$. The initial values for the estimates of the leader's velocity are all zero. The control parameter is chosen as $\gamma=0.04$ and $\Gamma_i=5I_3$, $i=1,\ldots,4$. Fig. 1 shows the trajectories of the leader and the followers. Clearly, all followers move cohesively with the leader without colliding with each other. Fig. 2 shows the velocity of the followers and the leader. It can be seen that the velocities of the followers converge to that of the leader and all agents move with the same velocity. There are two new edges added to the graph and no edge is lost.

\begin{figure}[hhhhtb]
\centering
\includegraphics[width=0.35\textwidth]{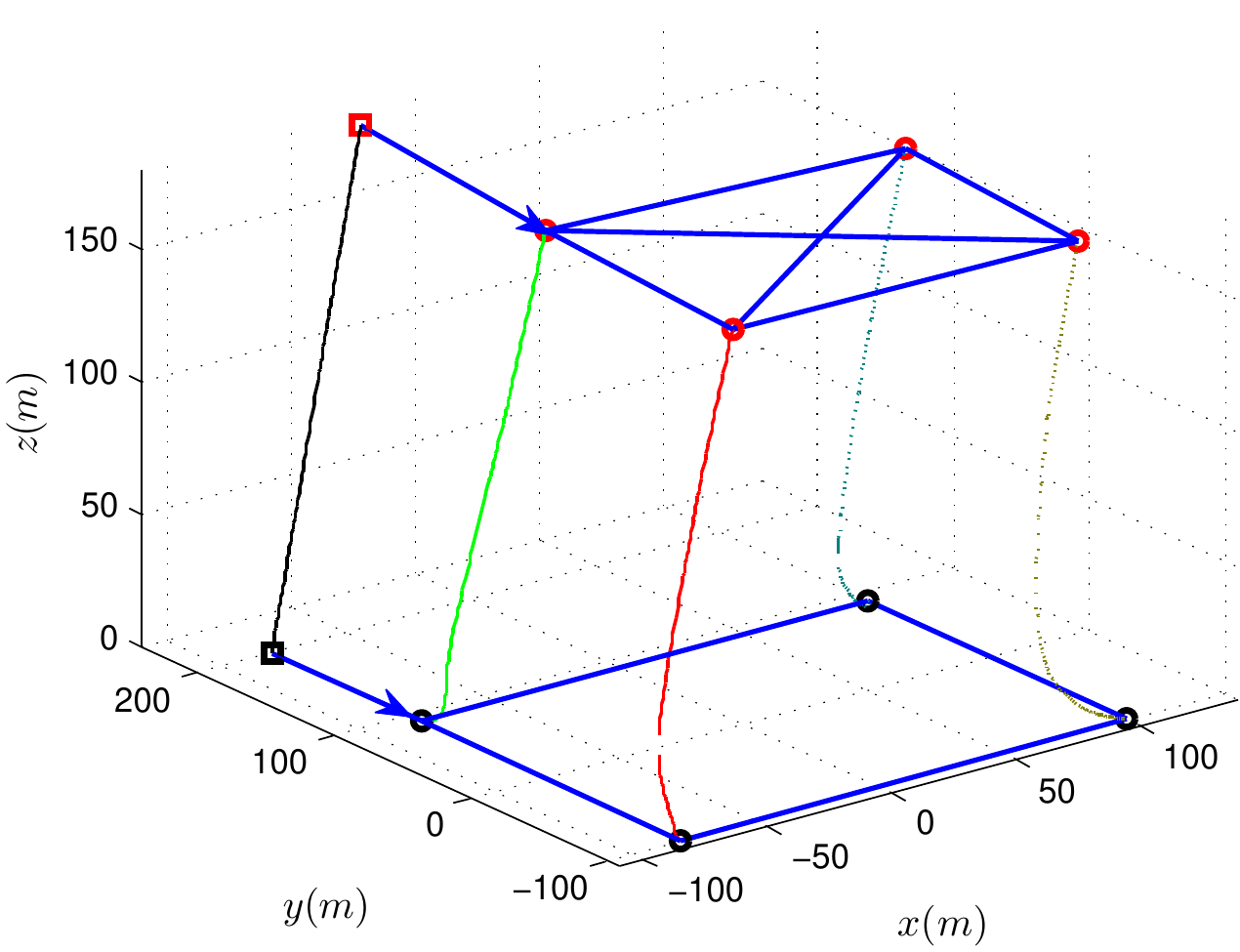}
\caption{The trajectories of the followers and the leader in the first case. The leader is represented as a square while the followers are represented as circles. An edge between two followers denotes that the two are neighbors, and an arrow from the leader to a follower denotes that the leader is a neighbor of the follower.}
\end{figure}

\begin{figure}[hhhhtb]
\centering
\includegraphics[width=0.35\textwidth]{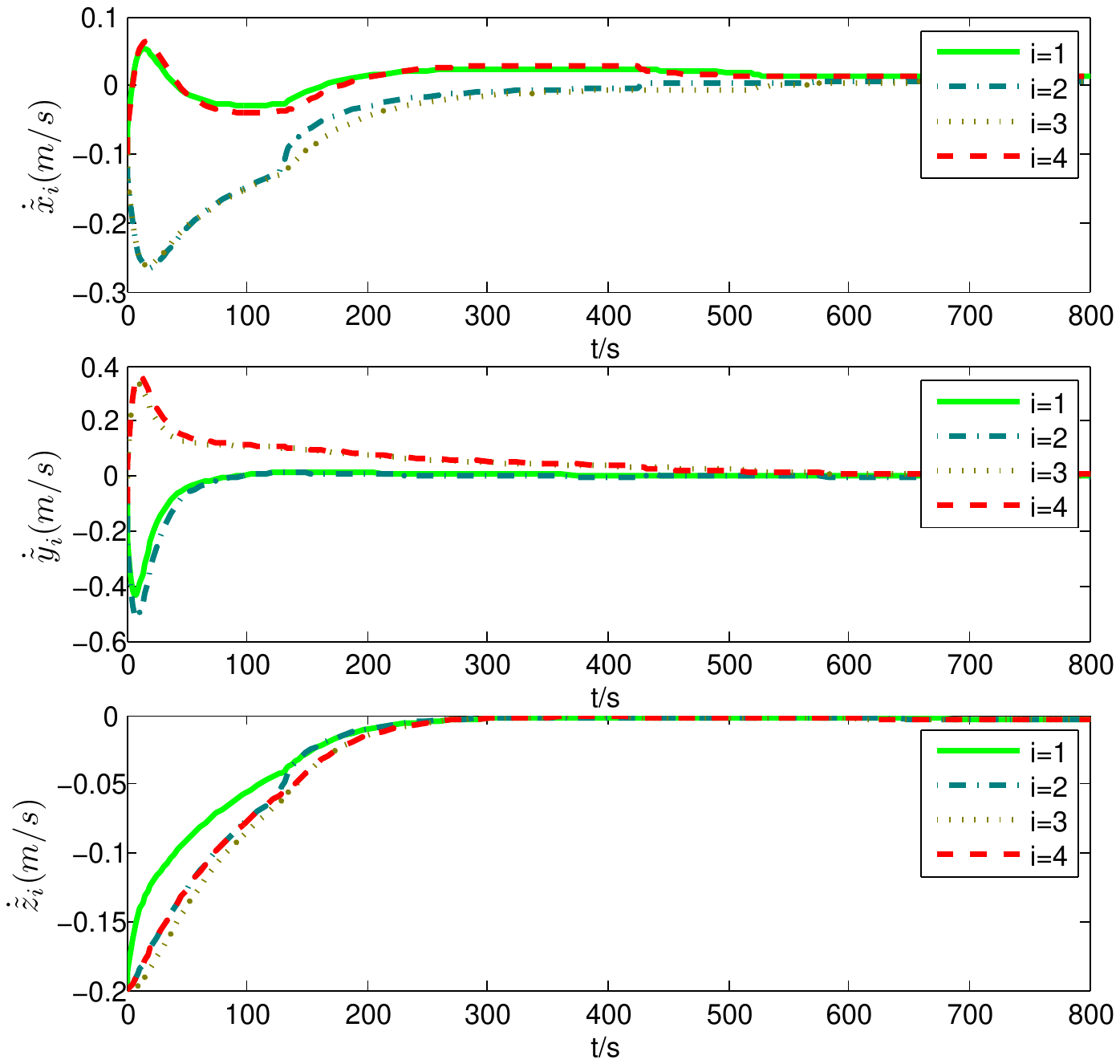}
\caption{The velocity errors between the followers and the leader using \eqref{input-cons-vel}-\eqref{par-estm-hat-cons-vel}.}
\end{figure}

In the second case, we simulate the case where the leader has a varying velocity under the control algorithm \eqref{input-varying-vel}-\eqref{par-estm-varying-vel}. The initial states of the followers are chosen as above and the leader's velocity is chosen as $\dot{q}_0(t)=[0.1\sin(\frac{2\pi}{60}t),0.1\cos(\frac{2\pi}{60}t),0.2]^T$. The initial position of the leader is chosen as $q_0(0)=[-80,200,0]^T$. The control parameters are chosen as $\alpha=0.04$, and $\Gamma_i=5 I_3$, $i=1,\ldots,4$. We use $\tanh(1000\cdot)$ to replace the function $\mbox{sgn}(\cdot)$.
Fig. 3 shows the trajectories of the followers and the leader. The agents maintain the initial connectivity while avoiding collisions. Fig. 4 shows that each follower eventually moves with the same velocity as the leader. Similarly, there are two new edges added to the graph and no edge is lost.

\begin{figure}[hhhhtb]
\centering
\includegraphics[width=0.35\textwidth]{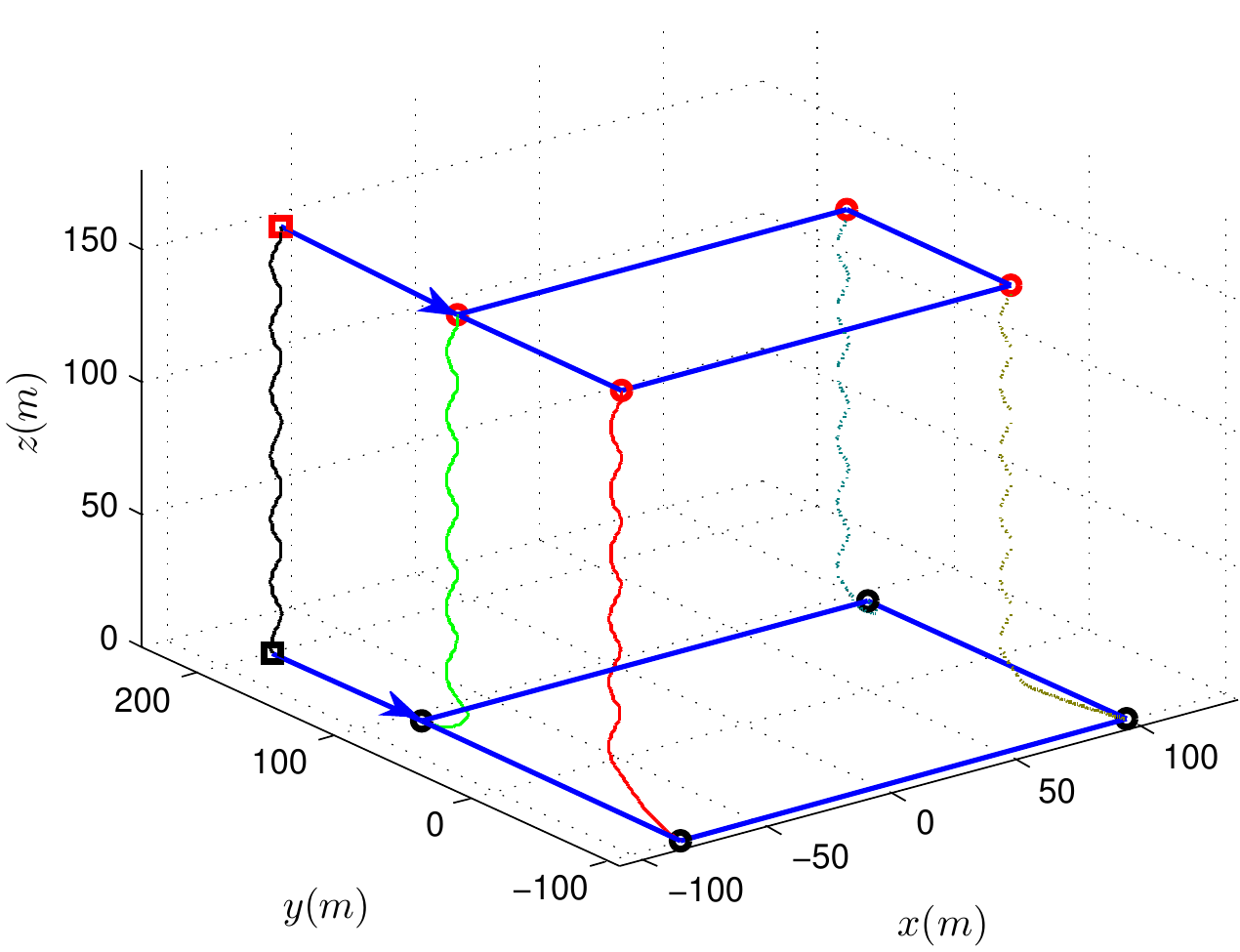}
\caption{The trajectories of the followers and the leader in the second case.}
\end{figure}

\begin{figure}[hhhhtb]
\centering
\includegraphics[width=0.35\textwidth]{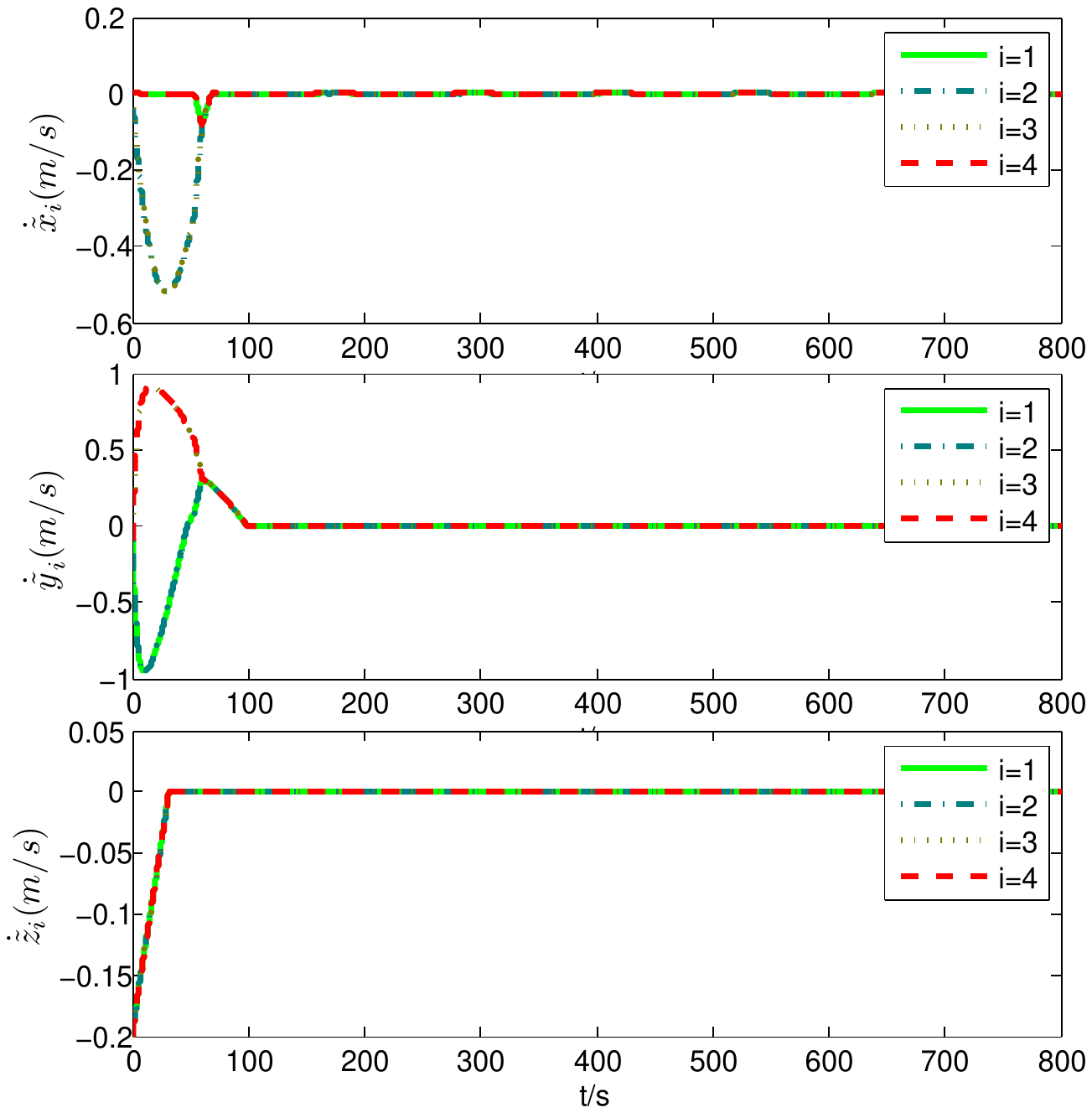}
\caption{The velocity errors between the followers and the leader using \eqref{input-varying-vel}-\eqref{par-estm-varying-vel}.}
\end{figure}

%

In the third case, we simulate the case where the leader has a varying velocity under the fully distributed control algorithm \eqref{equ28}-\eqref{equ34}. Here the initial states and the leader's trajectory are chosen as the second case. The control parameter is chosen as $\Gamma_i=5I_3$, $\gamma_{1i}=\gamma_{2i}=0.003$, $i=1,\ldots,4$.
Fig. 5 shows the trajectories while Fig. 6 shows the velocities of the followers and the leader. It can be seen that the leader-following flocking is achieved and there is no edge added or lost.

\begin{figure}[hhhhtb]
\centering
\includegraphics[width=0.35\textwidth]{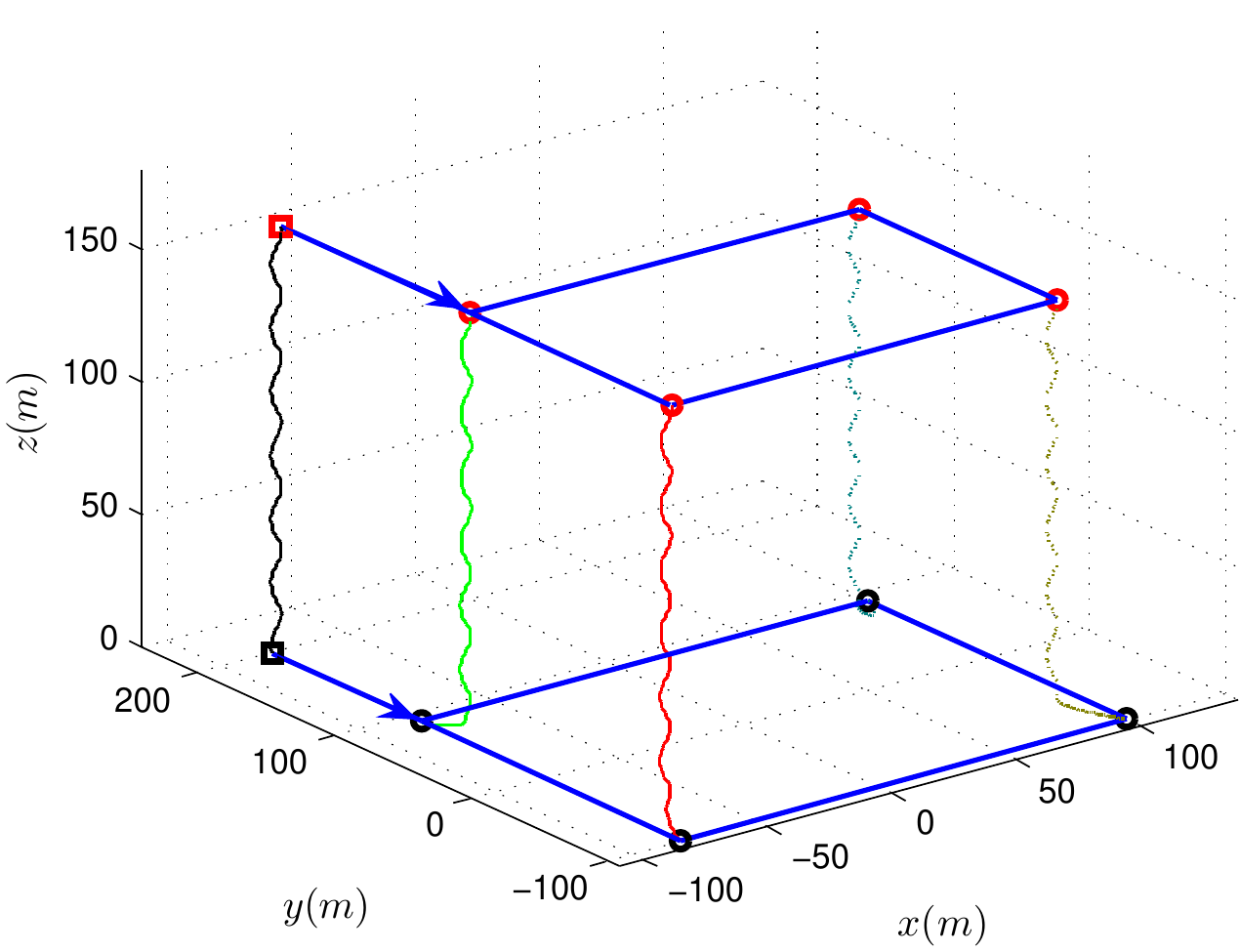}
\caption{The trajectories of the followers and the leader in the third case.}
\end{figure}

\begin{figure}[hhhhtb]
\centering
\includegraphics[width=0.35\textwidth]{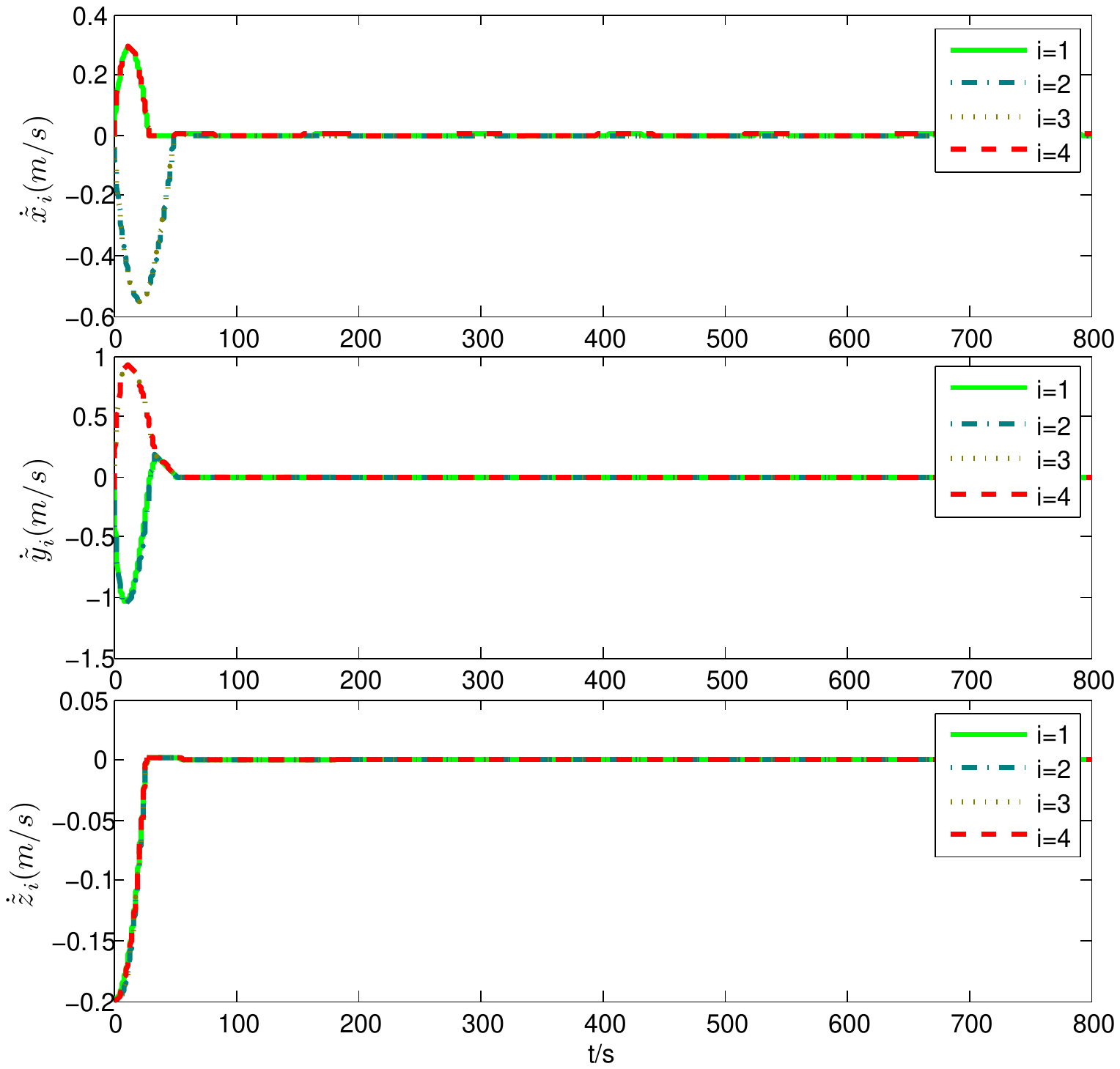}
\caption{The velocity errors between the followers and the leader using \eqref{equ28}-\eqref{equ34}.}
\end{figure}

%

\section{CONCLUSIONS}



In this paper, the distributed leader-follower flocking problem has been studied. The agents' models are described by Lagrange dynamics with unknown but constant parameters. Two cases for the leader have been considered: i) the leader has a constant velocity, and ii) the leader has a varying velocity. In both cases the leader is a neighbor of only a group of followers and the followers interact with only their neighbors defined by a proximity graph. In the second case we also relaxed the assumption of global information for parameter determination and proposed a fully distributed control algorithm. All proposed control algorithms require only one-hop neighbors' information and have been shown to achieve connectivity maintenance, collision avoidance, and velocity matching with a moving leader. Numerical simulations have also been presented to illustrate the theoretical results.

\bibliographystyle{elsart-num}    
\bibliography{paper-refs}           



\end{document}